\documentclass[12pt]{article}
\usepackage{amsfonts}
\input xy
\xyoption{all}
\newtheorem{theorem}{Theorem}[section]
\newtheorem{lemma}[theorem]{Lemma}

\newtheorem{proposition}[theorem]{Proposition}
\newtheorem{corollary}[theorem]{Corollary}

\newcommand{\Skew}{\mbox{\rm Skew}}

\newcommand{\mZ}{\mbox{$\mathbb Z$}}
\newcommand{\cA}{\mbox{$\mathcal A$}}
\newcommand{\cQ}{\mbox{$\mathcal Q$}}

\begin{document}
\title{Mutation classes of skew-symmetric $3 \times 3-$matrices} 
\author{Ibrahim Assem \and Martin Blais \and Thomas Br\"ustle \and  Audrey Samson\thanks{Ibrahim Assem is partially supported by NSERC of Canada. Martin Blais and Audrey Samson were working under a summer research fellowship of NSERC. Thomas Br\"ustle is partially supported by NSERC, by Bishop's University and the University of Sherbrooke.}}
  \date{}
\maketitle
\begin{center}\large Dedicated to the memory of A.V.Roiter
\end{center}
\bigskip

\begin{abstract}
In this paper, we establish a bijection between the set of mutation classes of mutation-cyclic skew-symmetric integral $3 \times 3-$matrices and the set of triples of integers $(a,b,c)$ such that $2 \le a \le b \le c$ and $ab \ge c$.
We also give an algorithm allowing to verify whether a matrix is mutation-cyclic or not. We prove that these two cases are not intertwined.
\end{abstract}

\section{Introduction}
Cluster algebras have been introduced  and studied by Fomin and Zelevinsky in \cite{FZ1,FZ2,BFZ}.
In particular, it was shown in  \cite{FZ2}  that every cluster algebra of finite type is acyclic, and corresponds to a Dynkin quiver. In \cite{BMRRT}, a fruitful connection between acyclic cluster algebras and representations of quivers has been established. In general, there is no known distinction between cyclic and acyclic cluster algebras. 
The objective of the present note is to study the first non-trivial situation, that of the (coefficient-free) cluster algebras of rank three which are given by a square skew-symmetric matrix.
\medskip

In general, for any square skew-symmetric integral matrix $B$, we denote by $\cA(B)$ the associated coefficient-free cluster algebra, as in \cite{FZ1}. 
We say that $\cA(B)$ has rank $n$ if $B$ is an $n \times n-$matrix.
The algebra $\cA(B)$ is constructed using mutations on $B$, thus depends not on $B$ itself, but rather on its mutation class within the set $\Skew_n(\mZ)$ of the skew-symmetric integral $n \times n-$matrices.
Every matrix $B=[b_{ij}]$ in $\Skew_n(\mZ)$ determines a quiver $Q_B$ having $\{1, \ldots , n\}$ as set of points, and $b_{ij}$ arrows from $i$ to $j$ whenever $b_{ij} > 0$. 
Thus, mutations on $B$ can equivalently be expressed as mutations on the quiver $Q_B$ (see (2.1)).
The cluster algebra ${\cA}(B)$ is called {\em acyclic} if there is a matrix in the mutation class of $B$ whose quiver  is acyclic, and otherwise it  is called cyclic (see  \cite{FZ2}). 
We say that a  matrix $B \in \Skew_n(\mZ)$ is {\em cyclic} (or {\em acyclic}) if the quiver $Q_B$ is so, and we say that $B$ is {\em mutation-cyclic} (or {\em mutation-acyclic}) if the corresponding cluster-algebra  $\cA(B)$ is cyclic (or acyclic, respectively).
\medskip

In this paper, we consider the case where $n=3$: this corresponds to quivers with three points. As we see in (2.2) below, the mutation class of a connected quiver with three points always contains a cyclic representative which is determined by three positive integral parameters $a,b,c$ corresponding to the number of arrows on each side.
Up to orientation, we may suppose that $a \le b \le c$.
Accordingly, the cyclic matrices $B \in \Skew_3(\mZ)$ are determined, up to transposition and simultaneous permutation of rows and columns, by the three parameters  $a \le b \le c$.
Our first theorem characterizes mutation-cyclic quivers (or, equivalently, matrices), in terms of these parameters.
Note that we consider all matrices up to transposition and permutation.
 
\begin{theorem}\label{thmA}
There exists a bijection between the set of mutation classes of  mutation-cyclic matrices  in $\Skew_3(\mZ)$ and the set of  triples $(a,b,c)$ of integers such that $2 \le a \le b \le c$ and $ab \ge 2c$. \end{theorem}

Our proof gives at the same time a handy algorithm allowing to verify whether a given matrix (or the corresponding quiver) is mutation-cyclic or not.
In our second theorem, we show that the two cases (mutation-cyclic and mutation-acyclic) are not intertwined.

\begin{theorem}\label{thmB}
Let $2 \le a \le b \le c$. Then there exists a unique  integer $c_0 \in [ab-b,ab-1]$ such that a cyclic matrix  in $\Skew_3(\mZ)$ which is represented by the triple $(a,b,c)$  is mutation-acyclic if and only if $c > c_0$.
\end{theorem}

The proofs of these theorems are purely combinatorial. 
The problem of characterizing the mutation-cyclic  skew-symmetric $3 \times 3 -$matrices is also considered in \cite{BBH} from the geometrical point of view.
In \cite{K}, the wild cluster-tilted algebras having three isomorphism classes of simple modules are studied. The quivers of those algebras yield the mutation-acyclic skew-symmetric $3 \times 3 -$matrices which are mutation-infinite.

The paper is organized as follows. In section 2, after a brief preliminary discussion, we prove Theorem \ref{thmA} and state our algorithm. Section 3 contains the proof of Theorem \ref{thmB} and ends with some examples.

%===============================================================================
\section{Mutation classes of quivers with three points}
\subsection{Preliminaries}
Since our intuition is graphical, we work with quivers rather than with matrices. 
For any quiver $Q$, we denote by $|Q_1|$ the number of arrows of $Q$.
The map $B \mapsto Q_B$ yields a bijection between $\Skew_n(\mZ)$ and the set $\cQ_n$ of quivers with $n$ points having neither loops nor cycles of length two. 
We begin by discussing the effect of certain matrix operations on the corresponding quivers.

When constructing the algebra $\cA(B)$, the matrix $B$ is given only up to simultaneous row and column permutations. This corresponds to considering the quivers $Q_B$ up to isomorphism. 
The transposition $B \mapsto B^T$ of the matrix $B$ corresponds to forming the opposite quiver.
To formulate our results in a more concise way, we always consider quivers up to isomorphism and change of orientation.
The mutation at $k$ of a matrix $B$ as defined in \cite{FZ1} yields a mutation of a quiver $Q \in \cQ_n$ at the point $k$ as follows:

\begin{itemize}
\item[(1)] All arrows passing through $k$ are reversed.
\item[(2)] If $Q$ has $r_{ij}$ paths of length two from $i$ to $j$ passing through $k$, then we add $r_{ij}$ arrows from $i$ to $j$.
\item[(3)] We delete all pairs of arrows which form cycles of length two.
\end{itemize}

Every mutation is an involution on $\cQ_n$. Taking the reflexive and transitive closure yields an equivalence relation denoted by $\sim$.
If $Q'$ is obtained from $Q$ by a mutation at the point $k$, we also write $Q \stackrel{k}{\sim} Q'$.

%--2.2--------------------------------------------------------------------------------------------------
\subsection{}
When we deal with a cyclic quiver $Q \in \cQ_3$, we always represent it as 

$$ 
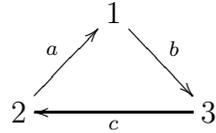
\xymatrix{ & 1 \ar^{b}[dr] & \\ 2  \ar^{a}[ur]  & & 3 \ar^{c}[ll]  }  $$
\medskip

where $a,b,c$ (called the {\em parameters} of $Q$) represent the number of arrows in the shown direction.
Up to isomorphism of quivers and replacing $Q$ by its opposite quiver, we can suppose that $0 < a \le b\le c$.

\begin{lemma}\label{lem-cyclic}
Let $Q$  be a connected quiver in $\cQ_3$, then
\begin{itemize}
\item[(1)] $Q$ is mutation-equivalent to a cyclic quiver.
\item[(2)] $Q$ is mutation-acyclic if and only if $Q$ is mutation-equivalent to the quiver $Q'$
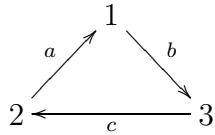
$$ \xymatrix{ & 1 \ar^{b}[dr] & \\ 2  \ar^{a}[ur]  & & 3 \ar^{c}[ll]  }  $$
where $0 < a \le b\le c$ and $ab \le c$.
\end{itemize}
\end{lemma}

{\bf Proof.}
(1) If $Q$ is acyclic, it is of the form 
$$ 
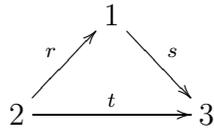
\xymatrix{ & 1 \ar^{s}[dr] & \\ 2  \ar^{t}[rr] \ar^{r}[ur]  & & 3  }  $$
where $r,s,t$ represent the numbers of arrows in the shown direction.
Since $Q$ is connected, we may assume that at most one of $r,s,t$ is zero. 
If $r$ and $s$ are non-zero, then mutation at $1$ yields the quiver $Q''$
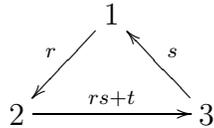
$$ \xymatrix{ & 1 \ar_{r}[dl] & \\ 2  \ar^{rs+t}[rr] & & 3  \ar_{s}[ul]  }  $$
which is cyclic because $rs+t \ge rs > 0$.

If $r$ and $t$ are non-zero, then mutation at $3$ yields the quiver $Q'''$
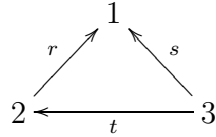
$$ \xymatrix{ & 1 & \\ 2 \ar^{r}[ur]& & 3  \ar_{s}[ul] \ar^{t}[ll]  }  $$
and thus, as before, mutation at $2$ yields a cyclic quiver. The case where $s$ and $t$ are non-zero is dual.
\medskip

(2) The quiver $Q'$ from (2) is clearly mutation-acyclic, since mutating at $1$ yields an acyclic quiver when $ab \le c$.
Conversely, starting with an acyclic quiver $Q$, we obtain in the proof of (1) the quivers $Q''$ and $Q'''$.
Then $Q''$ satisfies the required condition, as seen by setting $a=s, b=r$ and $c=rs+t$. Similarly for $Q'''$.
\hfill $\Box$ \bigskip

%--2.3--------------------------------------------------------------------------------------------------
\subsection{Mutation-finite quivers}
A quiver in $\cQ_n$ is called {\em mutation-finite} if its mutation class is a finite set.
The results in this subsection can also be derived from \cite{Seven} and \cite{SP}.

\begin{lemma}\label{lem-finite}
Let $Q$  be a connected quiver in $\cQ_3$, then
$Q$ is mutation-finite if and only if any cyclic quiver $Q'$ in its mutation class has at most two parallel arrows.
\end{lemma}
{\bf Proof.}
Sufficiency. There are only four connected cyclic quivers in $\cQ_3$ with at most two parallel arrows:
$$     
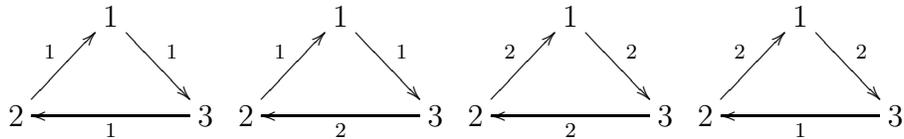
\xymatrix{ & 1 \ar^{1}[dr] & \\ 2  \ar^{1}[ur]  & & 3 \ar^{1}[ll]  } \;
      \xymatrix{ & 1 \ar^{1}[dr] & \\ 2  \ar^{1}[ur]  & & 3 \ar^{2}[ll]  } \;
       \xymatrix{ & 1 \ar^{2}[dr] & \\ 2  \ar^{2}[ur]  & & 3 \ar^{2}[ll]  } \;
        \xymatrix{ & 1 \ar^{2}[dr] & \\ 2  \ar^{2}[ur]  & & 3 \ar^{1}[ll]  }
$$
By \cite{SP}, the first three quivers are 4-bounded, hence mutation-finite. 
The last one does not satisfy the hypothesis, because mutation at 1 yields a cyclic quiver with 3 arrows from 3 to 2.

Necessity. Assume $Q' \sim Q$ with $Q'$ of the form
$$ 
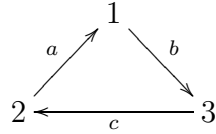
\xymatrix{ & 1 \ar^{b}[dr] & \\ 2  \ar^{a}[ur]  & & 3 \ar^{c}[ll]  }  $$
where $0 < a \le b \le c$ and $c \ge 3$. Mutating at 3 yields a quiver $Q''$ having $bc-a$ arrows from 1 to 2. Since $c \ge 3$ and $ b \ge a$, then $bc-a \ge 2a > a$. Thus $|Q_1''| > |Q_1'|$. 
Induction shows that we may repeat this procedure infinitely many times leading to quivers with more and more arrows. Therefore $Q$ is not mutation-finite.
\hfill $\Box$ \bigskip

\begin{corollary}\label{cor2.4}
There are exactly three mutation classes of mutation-finite connected quivers in $\cQ_3$, given by the following representatives:
$$  
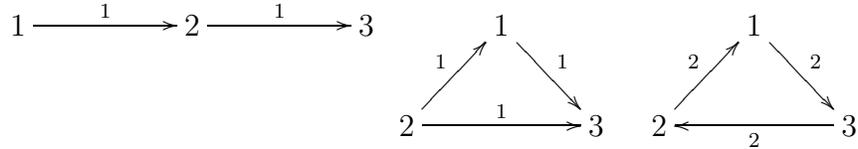
\xymatrix{ 1 \ar^{1}[rr] & &  2  \ar^{1}[rr]  & &  3  } \;
       \xymatrix{ & 1 \ar^{1}[dr] & \\ 2  \ar^{1}[ur]  \ar^{1}[rr]  & & 3   } \quad
       \xymatrix{ & 1 \ar^{2}[dr] & \\ 2  \ar^{2}[ur]  & & 3 \ar^{2}[ll]  } 
       $$
       Of these, only the last quiver is mutation-cyclic. Moreover, it is the only representative in its mutation class.
\end{corollary}

%--2.4--------------------------------------------------------------------------------------------------
\subsection{}
Consider a cyclic quiver given by parameters $0 < a \le b \le c$. Then  Lemma \ref{lem-cyclic}, (2) implies that the quiver is mutation-acyclic if $a=1$. 
Thus, since we are interested in mutation-cyclic quivers, we may assume that $(a,b,c) \in \mZ^3$ satisfy $2 \le a \le b \le c$. The following lemma shows when such a quiver stays cyclic under one mutation. It also describes how the parameters change.
\medskip

\begin{lemma}\label{lem-mut}
Let $Q$ be the cyclic quiver
$$ 
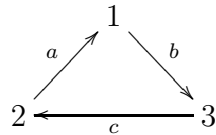
\xymatrix{ & 1 \ar^{b}[dr] & \\ 2  \ar^{a}[ur]  & & 3 \ar^{c}[ll]  }  $$
with  $2 \le a \le b \le c$, and $Q'$ be obtained from $Q$ by a mutation at a point $k \in \{1,2,3\}$.
\begin{itemize}
\item[(1)] If $k=1$, then $Q'$ is cyclic with $|Q_1'| \ge |Q_1|$ if and only if $ab \ge 2c$.
\item[(2)] If $k=2$, then $Q'$ is cyclic with the parameters satisfying 
\newline $2 \le a \le c \le ac-b$.
\item[(3)] If $k=3$, then $Q'$ is cyclic with the parameters satisfying 
\newline $2 \le b \le c \le bc-a$.
\end{itemize}
\end{lemma}
{\bf Proof.}
Assume first $k=1$.
\begin{itemize}
\item[i)] If $c \ge ab$, then $|Q_1'| = a+b+(c-ab) = |Q_1| - ab < |Q_1|.$
\item[ii)] If $c < ab$, then $|Q_1'| = a+b+(ab-c) = |Q_1|+(ab-2c).$
\end{itemize}
Thus, $Q'$ is cyclic when $ab >c$, and $|Q_1'| \ge |Q_1|$ if and only if $ab \ge 2c$.
\medskip

Assume now $k=2$. Then the numbers of arrows of $Q'$ are $a,c$ and $ac-b$, respectively. 
Clearly $2 \le a \le c$. 
Moreover, since $a \ge 2$ and $c\ge b$ we have  $ac \ge c+c \ge b+c$, thus $ac-b \ge c$ which is the inequality we wanted.
The proof is similar for $k=3$.
\hfill $\Box$ \bigskip

{\bf Remark.} Parts (2) and (3) of the preceding lemma show that mutation in a point opposite to one of the parameters $a$ or $b$ necessarily yields a new maximal parameter. This new parameter is strictly greater than the two other parameters unless  $a=2$ and $b=c$ in (2) or  $a=b=c=2$ in case (3).

\begin{proposition}\label{prop1}
Let $Q$ be the cyclic quiver
$$ 
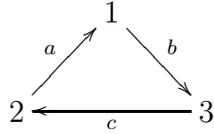
\xymatrix{ & 1 \ar^{b}[dr] & \\ 2  \ar^{a}[ur]  & & 3 \ar^{c}[ll]  }  $$
with  $2 \le a \le b \le c$.
 If $ab \ge 2c$, then $Q$ is mutation-cyclic and, for any $Q' \sim Q$, we have $|Q_1'| \ge |Q_1|$. 
\end{proposition}

{\bf Proof.}
If $Q' \sim Q$, there exists  a sequence of mutations 
$$Q = Q^{(0)} \stackrel{k_1}{\sim} Q^{(1)} \stackrel{k_2}{\sim}  \cdots \stackrel{k_n}{\sim} Q^{(n)} = Q' $$ 
Without loss of generality, we may assume this sequence to be reduced, that is, two consecutive mutations in this sequence are not inverse to each other (thus, $k_i \neq k_{i+1}$ for all $i$ ).
We show by induction on $n$ that the quiver $Q^{(n)}$ is cyclic (and so are all quivers $Q^{(i)}$ in the sequence) and that  $|Q^{(n)}| \ge |Q^{(n-1)}|$.

Lemma \ref{lem-mut} shows that $Q^{(1)}$ is cyclic with $|Q_1^{(1)}| \ge |Q_1|$, thus we consider  the induction step.
Suppose $n \ge 2$, and let $Q^{(n-1)}$ be the quiver
$$ 
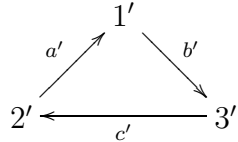
\xymatrix{ & 1' \ar^{b'}[dr] & \\ 2'  \ar^{a'}[ur]  & & 3' \ar^{c'}[ll]  }  $$
which   is cyclic by induction. Up to duality, we can suppose  $ a' \le b' \le c'$, and by induction, the numbers of arrows cannot decrease by mutations, thus $2 \le a'$.
The quiver $Q^{(n)}$ is obtained from $Q^{(n-1)}$ by mutation at $k_n$.
If $k_n \neq 1'$ then the statement follows by lemma \ref{lem-mut}.
Assume now that $k_n=1'$. 
Since the sequence above  is supposed to be minimal, the mutation preceding $k_n$ satisfies $k_{n-1} \neq 1'$.
This means that one of the values $a'$ or $b'$ has been changed when going from $Q^{(n-2)}$ to $Q^{(n-1)}$.
But we know that $ a' \le b' \le c'$, therefore the remark following  lemma \ref{lem-mut} shows that we are in one of the cases  $a'=2$ and $b'=c'$ or   $a'=b'=c'=2$.
If $a'=b'=c'=2$, then all the quivers $Q^{(i)}$  are isomorphic, and the statement holds. 
If  $a'=2$ and $b'=c'$, then mutation at $1'$ transforms the number $c'$ in $Q^{(n-1)}$ into $a'b'-c'=2b'-b'=b'=c'$, thus it stays the same, which implies that $Q^{(n)}$  is cyclic and $|Q_1'| \ge |Q_1|$, which we wanted to show.
\hfill $\Box$ \bigskip

%--2.5--------------------------------------------------------------------------------------------------
\subsection{}

In the sequel, we call {\em root} (of its mutation class) a cyclic quiver $Q \in \cQ_3$ whose parameters $(a,b,c)$ satisfy $2 \le a \le b\le c$ and $ab \ge 2c$.
\medskip

The previous results yield the following algorithm which decides  whether a given connected quiver $Q \in \cQ_3$ is mutation-cyclic or not:

\begin{itemize}
\item[(1)] If $Q$ is acyclic, stop.
\item[(2)] Otherwise, $Q$ is a cyclic quiver with ordered parameters $0 \le a \le b \le c$.
If $ab \le c$, then $Q$ is mutation-acyclic by lemma \ref{lem-cyclic}, stop.
\item[(3)] Perform a mutation at the point opposite to $c$. If the number of arrows has decreased, go back to step (1). 
Otherwise, $Q$ is a root by lemma \ref{lem-mut}, and thus $Q$ is mutation-cyclic.
\end{itemize}
\medskip

This procedure must clearly stop after finitely many steps, since we deal with strictly decreasing sequences of natural numbers. 
\medskip

{\bf Proof of theorem \ref{thmA}.}
We have shown in proposition \ref{prop1} that every root is mutation-cyclic.
Conversely, let $Q$ be a mutation-cyclic quiver. Applying the algorithm above, we find a root in the mutation class of $Q$. Up to duality, this root is uniquely described by its parameters $a \le b \le c$.
\hfill $\Box$ \bigskip
\medskip

%===============================================================================
\section{Separating the cyclic case from the acyclic}
\begin{lemma}\label{lem-mutcyclic}
The quiver $Q$
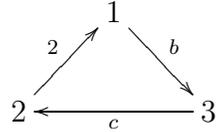
$$ \xymatrix{ & 1 \ar^{b}[dr] & \\ 2  \ar^{2}[ur]  & & 3 \ar^{c}[ll]  }  $$
with $2 \le b \le c$ is mutation-cyclic if and only if $b=c$.
\end{lemma}

{\bf Proof.}
Sufficiency follows from proposition \ref{prop1}, so we only show necessity.
Asssume that $Q$ is mutation-cyclic. 
By applying a sequence of mutations at the point opposed to the maximum, we reach a root $Q'$ 
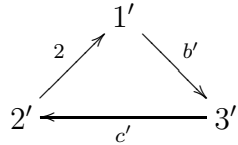
$$ \xymatrix{ & 1' \ar^{b'}[dr] & \\ 2'  \ar^{2}[ur]  & & 3' \ar^{c'}[ll]  }  $$
with $2 \le b' \le c'$ and $2b' \ge 2c'$ (because $Q'$ is a root). Thus $b' = c'$.
But mutations at $1'$ or $2'$ do not change the quiver $Q'$. This implies that $Q'=Q$, and so $b=b'=c'=c$.
\hfill $\Box$ \bigskip

%----------------------------------------------------------------------------------------------------
\begin{lemma}\label{lem3.2}
Let $Q$ be the quiver
$$ 
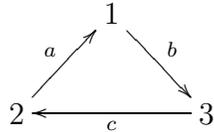
\xymatrix{ & 1 \ar^{b}[dr] & \\ 2  \ar^{a}[ur]  & & 3 \ar^{c}[ll]  }  $$
with $2\le a \le b \le c$.
\begin{itemize}
\item[(1)] If $c \le ab-b$, then $Q$ is mutation-cyclic.
\item[(2)]  If $c \ge ab-1$, then $Q$ is mutation-acyclic.
\end{itemize}
\end{lemma}

{\bf Proof.}
Mutating at $1$ yields a quiver $Q'$ with $ab-c$ arrows between $2$ and $3$.
\begin{itemize}
\item[(1)] If $ab-c \ge b$ (or, equivalently, $c \le ab-b$), then the new maximal number is $ab-c$. Since the maximum did not change its position, it follows from lemma \ref{lem-mut} that $Q$ or $Q'$ is a root. In particular, $Q$ is mutation-cyclic.
\item[(2)]  If $ab-c \le 1$ (or, equivalently, $c \ge ab-1$), then $Q$ is mutation-acyclic by lemma 
\ref{lem-cyclic}, (2).
\end{itemize}
\hfill $\Box$ \bigskip

%----------------------------------------------------------------------------------------------------
\begin{lemma}\label{lem3.3}
Consider the homogeneous difference equation
$$ S_{n+2} = a S_{n+1} - S_n$$
(with $a \ge 2, S_0=0$ and $S_1=1$), then
\begin{itemize}
\item[(1)]  $S_n = \left\{ \begin{array}{lr}
\frac{1}{2^n\sqrt{a^2-4}} \left[ (a+\sqrt{a^2-4})^n - (a-\sqrt{a^2-4})^n\right] & \mbox{ if } a\ge 3 \\
n  & \mbox{ if } a =2
\end{array}  \right. $
\item[(2)]  The sequence $(S_n)_{n \ge 0}$ is strictly increasing
\item[(3)]  $S_n \ge a$ for any $n \ge 2$.
\end{itemize}
\end{lemma}

{\bf Proof.}
\begin{itemize}
\item[(1)] This is straightforward and left to the reader.
\item[(2)]  Use induction on $n$: clearly $S_0 < S_1$. Assume $S_n < S_{n+1}$, then $S_{n+2} = aS_{n+1} - S_n > aS_{n+1} - S_{n+1} \ge S_{n+1}$ because $a \ge 2$.
\item[(3)] Since $S_2=a$, this follows from $(2)$.\hfill $\Box$ \bigskip
\end{itemize}

%---------------------------------------------------------------------------------------------------
{\bf Proof of theorem \ref{thmB}.}
Let $Q$ be the quiver
$$ 
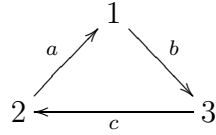
\xymatrix{ & 1 \ar^{b}[dr] & \\ 2  \ar^{a}[ur]  & & 3 \ar^{c}[ll]  }  $$
with $2\le a \le b \le c$.
By lemma \ref{lem-mutcyclic}, the statement holds if $a=2$ (with $c_0 = b$).
 Assume $a \ge 3$ and that such a $c_0$ does not exist.
Then there exist quivers $Q,Q'$ 
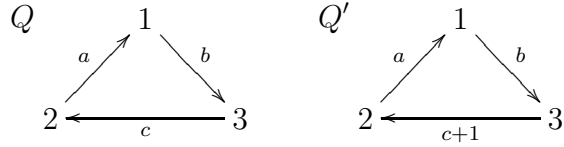
$$Q \xymatrix{ & 1 \ar^{b}[dr] & \\ 2  \ar^{a}[ur]  & & 3 \ar^{c}[ll]  }  \qquad
Q' \xymatrix{ & 1 \ar^{b}[dr] & \\ 2  \ar^{a}[ur]  & & 3 \ar^{c+1}[ll]  }  $$
with $3 \le a \le b \le c$ such that $Q$ is mutation-acyclic and $Q'$ is mutation-cyclic. We show by induction that the quiver $Q$ never reaches an acyclic representative. 
This will yield a contradiction which implies our statement.

More precisely, we prove by induction on $n\ge 1$ that, after $n$ mutations at the point opposed to the maximum, we obtain respectively the quivers $Q^{(n)}$ and $Q'^{(n)}$
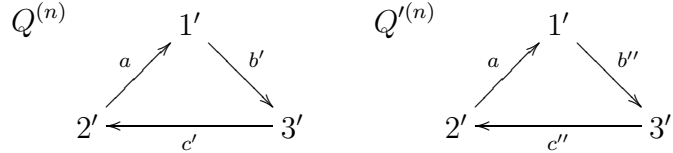
$$Q^{(n)} \xymatrix{ & 1' \ar^{b'}[dr] & \\ 2'  \ar^{a}[ur]  & & 3' \ar^{c'}[ll]  }  \qquad
Q'^{(n)} \xymatrix{ & 1' \ar^{b''}[dr] & \\ 2'  \ar^{a}[ur]  & & 3' \ar^{c''}[ll]  }  $$
with $S_n=b'-b'', S_{n-1}=c'-c'', \max(a,b',c')=c'$ and $ \max(a,b'',c'')=c''$ where $S_n$ is as in lemma \ref{lem3.3} and $Q^{(n)}$ is still a cyclic quiver (note that the parameter $a$ does not change).

Assume first $n=1$. After one mutation at the point opposed to the maximum, we obtain respectively the quivers
$$Q^{(1)} 
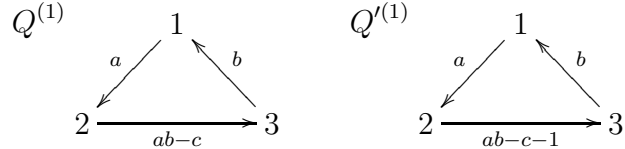
\xymatrix{ & 1 \ar_{a}[dl] & \\ 2  \ar_{ab-c}[rr] & & 3  \ar_{b}[ul]  }  \qquad
 Q'^{(1)} \xymatrix{ & 1 \ar_{a}[dl] & \\ 2  \ar_{ab-c-1}[rr] & & 3  \ar_{b}[ul]  } $$
If $Q^{(1)}$ is not cyclic, then $ab-c \le 0$, hence $ab-c-1 < 0$ and so $Q'^{(1)}$ is not cyclic either.
This is a contradiction to the hypothesis that $Q'$ is mutation-cyclic, hence $Q^{(1)}$ is a cycle.
In order to pursue the algorithm, we must have $b > ab-c$ (otherwise, the next mutation yields $Q^{(2)}=Q$ and so $Q$ or $Q^{(1)}$ is a root, contradicting the hypothesis that $Q$ is mutation-acyclic).
Therefore, $b > ab-c-1$.
Thus, in $Q^{(1)}$ and $Q'^{(1)}$, the maxima correspond to the same sides. 
Moreover, $S_1 = (ab-c)-(ab-c-1) =1$ and $S_0 = b-b=0$.
\medskip

Assume the statement holds for $n$. Mutating at $1'$ (which is opposed to the respective maxima $c'$ and $c''$) yields respectively the quivers 
$$Q^{(n+1)} \xymatrix{ & 1' \ar_{a}[dl] & \\ 2'  \ar_{ab'-c'}[rr] & & 3'  \ar_{b'}[ul]  }  \qquad
 Q'^{(n+1)} \xymatrix{ & 1' \ar_{a}[dl] & \\ 2'  \ar_{ab''-c''}[rr] & & 3'  \ar_{b'}[ul]  } $$
We first note that
$(ab'-c') - (ab''-c'') = a(b'-b'') - (c'-c'') = a S_n-S_{n-1}= S_{n+1}$.
Also, $Q^{(n+1)}$ is a cycle. Indeed, if this is not the case, then $ab'-c' \le 0$. Since $S_{n+1}\ge 0$ (by lemma \ref{lem3.3}), we have $ab''-c'' \le 0$, contradicting the hypothesis that $Q'$ is mutation-cyclic.
\medskip

We now determine  the maximal parameter in $Q^{(n+1)}$.

1) Assume $\max (a,b',ab'-c') =a$. By lemma \ref{lem3.3}, we have
$$(ab'-c')-(ab''-c'')=S_{n+1} \ge a$$
hence $ab''-c'' \le (ab'-c') -a \le 0$, contradicting the hypotesis that $Q'$ is mutation-cyclic.

2)  Assume $\max (a,b',ab'-c') =ab'-c'$. Then, mutating at the point opposed to the maximum yields  
 $Q^{(n+2)}=Q^{(n)},$ thus $Q$ is mutation-cyclic, a contradiction.
 
 Therefore $\max(a,b',ab'-c') =b'$. We claim that also $\max(a,b'',ab''-c'') =b''$. 
 Assume first that $a > b''$. Since $ab''-c'' \ge 2$, we have
 $$\begin{array}{lcl}
 a(ab'-c')=a[(ab''-c'')+S_{n+1}] & = & a(ab''-c'')+aS_{n+1}\\
  & \ge & 2a+ a S_{n+1}\\
  & > & 2b''+a S_n\\
   & > & 2b''+2S_n = 2b'.
  \end{array}$$
On the other hand, $S_{n+1} = (ab'-c')-(ab''-c'') \ge 0$ yields $ab'-c' \ge ab''-c'' \ge 2$.
These inequalities show, by Theorem \ref{thmA}, that $Q^{(n+1)}$ is a root. Thus $Q$ is mutation-cyclic, a contradiction. Therefore $a \le b''$.

Since $b' > ab'-c'$, we have 
$b''=b'-S_{n} > (ab'-c')-S_{n+1}=ab''-c''$.

This shows that $b''=\max(a,b'',ab''-c'')$, completes the proof of the induction statements and thus establishes the theorem.
\hfill $\Box$ \bigskip
\medskip

%----------------------------------------------------------------------------------------------------
{\bf Example 1.}
Let $Q$ be the quiver 
$$ 
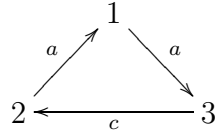
\xymatrix{ & 1 \ar^{a}[dr] & \\ 2  \ar^{a}[ur]  & & 3 \ar^{c}[ll]  }  $$
with $2\le a \le  c$.
Then we claim that $c_0=a^2-2.$

Indeed, applying a mutation at $1$ yields the quiver $Q'$ 
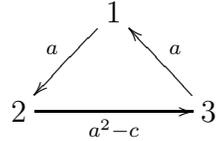
$$  \xymatrix{ & 1 \ar_{a}[dl] & \\ 2  \ar_{a^2-c}[rr] & & 3  \ar_{a}[ul]  } $$

1) Suppose $c \le a^2-2$. If $a^2-c$ is the maximum, then the next mutation at the point opposed to the 
maximum yields $Q$ so that $Q$ or $Q'$ is a root and $Q$ is mutation-cyclic. 
If not, then $a$ is the maximum and moreover $a(a^2-c) \ge 2a$ so that $Q'$ is a root, whence $Q$ is again mutation-cyclic.

2) Suppose $c > a^2-2$. Then $a^2-c \le1$ and so $Q$ is mutation-acyclic by lemma \ref{lem-cyclic}, (2).
\medskip

%----------------------------------------------------------------------------------------------------
{\bf Example 2.}
Let $Q$ be the quiver 
$$ 
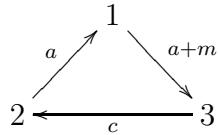
\xymatrix{ & 1 \ar^{a+m}[dr] & \\ 2  \ar^{a}[ur]  & & 3 \ar^{c}[ll]  }  $$
where $2 < a \le  c$ and $1 \le m \le 4$.
Then we claim that $c_0=a^2+am-3.$

Indeed, applying a mutation at $1$ yields the quiver $Q'$ 
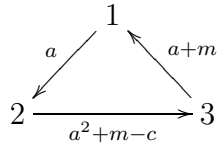
$$  \xymatrix{ & 1 \ar_{a}[dl] & \\ 2  \ar_{a^2+m-c}[rr] & & 3  \ar_{a+m}[ul]  } $$

1) Assume $c \le a^2+am-3$. If $a^2+am-c$ is the maximum, then the next mutation at the point opposed to the maximum yields $Q$ again, which is then mutation-cyclic.
If not, then $a+m$ is the maximum. 
If $a \ge 2m$, then $c \le a^2+am-3$ yields  $a(a^2+am-c) \ge 3a = 2a+a \ge 2 (a+m)$.
Thus, $Q'$ is a root and $Q$ is mutation-cyclic. On the other hand, if $a < 2m$, then we have
$$3 \le a < 2m \le 8$$
$$ a < a+m \le a+4$$
$$ a+m \le c \le a^2+am-c.$$
There are only finitely many quivers verifying these inequalities. A straight\-forward verification shows that in each case $Q$ is mutation-cyclic, as desired.

2) Assume $c > a^2+am-3$.
If $a^2+am-c \le 1$, then by lemma \ref{lem-cyclic}, (2) $Q$ is mutation-acyclic. Otherwise, $a^2+am-c=2$ and then the same conclusion follows from $(3.1)$.
 
\bigskip

{\sc D\'epartement de
  Math\'ematiques, 
  Universit\'e de Sherbrooke, Sherbrooke (Qu\'ebec), J1K 2R1, Canada}
 
 {\it  E-mail address: } {\tt ibrahim.assem@usherbrooke.ca}
 \medskip
 
{\sc D\'epartement de
  Math\'ematiques, 
  Universit\'e de Sherbrooke, Sherbrooke (Qu\'ebec), J1K 2R1, Canada}
 
 {\it  E-mail address: } {\tt Martin.V.Blais@usherbrooke.ca}
 \medskip
 
{\sc D\'epartement de
  Math\'ematiques, 
  Universit\'e de Sherbrooke, Sherbrooke (Qu\'ebec), J1K 2R1, Canada} \quad and 
  
  {\sc Department of Mathematics, Bishop's University,   2600 College St., Sherbrooke, Quebec, Canada J1M 0C8} 

 {\it  E-mail address: } {\tt thomas.brustle@usherbrooke.ca} \quad and 
 
 {\tt tbruestl@ubishops.ca}
 \medskip
 
{\sc D\'epartement de
  Math\'ematiques, 
  Universit\'e de Sherbrooke, Sherbrooke (Qu\'ebec), J1K 2R1, Canada}
 
 {\it  E-mail address: } {\tt Audrey.Samson@USherbrooke.ca}
\end{document}